\documentclass{amsart}
\usepackage{amssymb}
\usepackage{enumerate}
\usepackage{url}

\newcommand\N{\mathbb N}
\newcommand\Z{\mathbb Z}

\newcommand\R{\mathbb R}

\newcommand\ph\varphi
\newcommand\ps\psi
\newcommand\ep\varepsilon
\newcommand\rh\varrho
\newcommand\al\alpha
\newcommand\be\beta
\newcommand\ga\gamma
\newcommand\om\omega
\newcommand\ta\tau
\renewcommand\th\theta
\newcommand\de\delta
\newcommand\ze\zeta
\newcommand\ch\chi
\newcommand\et\eta
\newcommand\io\iota
\newcommand\la\lambda
\newcommand\si\sigma

\newcommand\Ga\Gamma
\newcommand\De\Delta
\newcommand\Th\Theta
\newcommand\La\Lambda
\newcommand\Si\Sigma
\newcommand\Ph\Phi
\newcommand\Ps\Psi
\newcommand\Om\Omega
\newcommand\A{\mathcal{A}}
\newcommand\OO{\mathcal{O}}
\newtheorem{theorem}{Theorem}[section]
\newtheorem{lemma}[theorem]{Lemma}

\newtheorem{corollary}[theorem]{Corollary}
\theoremstyle{definition}
\newtheorem{definition}[theorem]{Definition}

\newtheorem{example}[theorem]{Example}
\theoremstyle{remark}

\newcommand\M{\mathfrak{m}}
\DeclareMathOperator{\sos}{sos}

\begin{document}
\title[SOS APPROXIMATIONS VIA SIMPLE HIGH DEGREE PERTURBATIONS]
{SOS APPROXIMATIONS OF NONNEGATIVE POLYNOMIALS VIA SIMPLE HIGH DEGREE PERTURBATIONS}
\author{Jean B. Lasserre}
\address{LAAS-CNRS, 7 Avenue du Colonel Roche, 31077 Toulouse C\'edex 4, France. 
Also associate member of IMT, the Institute of Mathematics of Toulouse.}
\email{lasserre@laas.fr}

\author{Tim Netzer}
\address{Universit\"at Konstanz, Fachbereich Mathematik und Statistik, 78457 Konstanz, Germany}
\email{tim.netzer@uni-konstanz.de}

\keywords{Real algebraic geometry; positive polynomials; sum of squares; semidefinite programming; moment problem.}
\subjclass{12E05, 12Y05, 90C22, 44A60}
\date{\today}

\begin{abstract}
We show that every real polynomial $f$ nonnegative on $[-1,1]^{n}$ can be approximated in the $l_{1}$-norm of coefficients, by a sequence of polynomials $\{f_{\ep r}\}$ that are sums of squares. This complements the existence of s.o.s. approximations in the denseness result of Berg, Christensen and Ressel, as we provide a very simple and \textit{explicit} approximation sequence.

Then we show that if the Moment Problem holds for a basic closed
semi-algebraic set $K_S\subset\R^n$ with nonempty interior, then every polynomial
nonnegative on $K_S$ can be approximated in a similar fashion by
elements from the corresponding preordering.

Finally, we show that the degree of the 
perturbation in the approximating sequence depends on $\epsilon$ as well as
the degree and the size of coefficients of the nonnegative polynomial $f$, but {\it not}
on the specific values of its coefficients.
\end{abstract}

\maketitle

\section{Introduction}
\textit{Sums of squares} (s.o.s.) polynomials are not only of self-interest, but are also of primary importance for practical computation, especially in view of their numerous potential applications, notably in polynomial optimization; see e.g. \cite{las1,par,sch,schw}. Indeed, in the computational complexity terminology, checking whether a given polynomial is nonnegative is a NP-hard problem, whereas checking whether it is s.o.s. reduces to solving a (convex) \textit{semidefinite programming} (SDP) problem which (up to arbitrary precision) can be done in time polynomial in the input size of the problem; for more detail on semidefinite programming, the interested reader is refered to Vandenberghe and Boyd \cite{van}.

It has been known for some time that the cone of s.o.s.
polynomials is \textit{dense} (for the $l_{1}$-norm of
coefficients) in the cone of polynomials nonnegative on the unit
ball $[-1,1]^{n}\subset \R^{n}$; see e.g. Berg, Christensen and Ressel \cite{bcr}
and Berg \cite{ber}. However, \cite{bcr} is essentially an existence result.

\textit{Contribution.} Our contribution is threefold:

(i) We first provide an explicit and very
simple s.o.s. approximation of polynomials  nonnegative on the
unit ball $[-1,1]^{n}$. Namely, let
\begin{equation}\label{Th}\Th_{r}:=1+\sum_{j=1}^{n}X_{j}^{2r}\in\R[X_{1},\ldots,X_{n}].\end{equation}
Then, given $\ep>0$ and a polynomial $f\in\R[X_{1},\ldots,X_{n}]$
nonnegative on $[-1,1]^{n}$, the polynomial $f_{\ep
r}:=f+\ep\Th_{r}$ is s.o.s. provided $r$ is large enough, say $r\geq r(f,\epsilon)$. Of
course, $\Vert f_{\ep r}-f \Vert_{1}\rightarrow 0$ as $\ep
\rightarrow 0$. Although our result is not completely constructive (as $r(f,\epsilon)$ is not known),
it complements the pure existence result \cite{bcr}.

If $f$ is nonnegative on the ball $[-l,l]^{n}$ for some $l>0$,
then for every $\ep>0$, the polynomial
$f+\ep(1+\sum_{j=1}^{n}(X_{j}/l)^{2r})$ is s.o.s. provided $r$ is
sufficiently large (just use $x\mapsto g(x):=f(lx)\geq 0$ on
$[-1,1]^{n}$).

Note that the  representation $f+\ep\Th_{r}=q_{\ep r}$ for some
s.o.s. polynomial $q_{\ep r}$, is an obvious \textit{certificate}
of nonnegativity of $f$ on $[-1,1]^{n}$. Indeed, for
$x\in[-1,1]^{n}$, one has $$f(x)+\ep\Th_{r}(x)=q_{\ep r}(x)\geq
0,$$ provided $r$ is big enough. As all $\Th_{r}$ are bounded by
$n+1$ on $[-1,1]^{n}$, letting $\ep \downarrow 0$ yields $f(x)\geq
0$.

Our s.o.s. approximation result  states that to approximate
(uniformly on $[-1,1]^{n}$) a polynomial nonnegative on
$[-1,1]^{n}$, it is enough to slightly perturb by a small $\ep>0$
its (maybe zero) coefficients of some even power of marginal
monomials $\{X_{i}^{2r}\}$.

The method of the proof is quite  different and much simpler than
that of \cite{las3} for s.o.s. approximation of nonnegative
polynomials; in particular, it does \textit{not} use Nussbaum's
deep result on moment sequences \cite{nuss}. It also simplyfies
the approximating sequence obtained in \cite{net} in the spirit of
\cite{las3}.

In addition, if one fixes \textit{\`a priori}  the degree $r$ of
the perturbation $\Th_{r}$, we also characterize the minimum value
$\ep_{r}^{*}$ of the parameter $\ep$, to make $f+\ep\Th_{r}$ a
s.o.s. It is given by $$-\ep_{r}^{*}:=\min_{L}\left\{ L(f)\mid
L:\mathcal{A}_{2r}\rightarrow \R \mbox{ linear, }L(\Th_{r})\leq 1,
\ L(h^{2})\geq 0\  \forall h\in\mathcal{A}_{r} \right\},$$ where
$\mathcal{A}_{r}$ is the finite dimensional vector space of
polynomials of degree at most $r$.

(ii) We next obtain a similar approximation result for polynomials
nonnegative on certain semi-algebraic sets. For a finite set $S\subset\R[X_1,\ldots,X_n]$
of polynomials, denote by $K_{S}$ the associated basic closed
semi-algebraic set in $\R^{n}$, and by $T_{S}$ the preordering
generated by $S$. Assume that $K_{S}$ has nonempty interior
and the so called \textit{Moment Problem} holds for $S$, that is,
every linear form on $\R[X_{1},\ldots,X_{n}]$ which is
nonnegative on $T_{S}$ comes from a measure on $K_{S}$. Then every
polynomial $f$ nonnegative on $K_{S}$ is approximated in the $l_1$-norm by the same
sequence $\{f_{\ep r}\}$, which now lies in $T_{S}$. In addition, if one uses the perturbation
\begin{equation}\label{smalltheta}\th_{r}:=\sum_{i=1}^n\sum_{k=0}^r
\frac{X_{j}^{2k}}{k{\rm !}}\quad\in\R[X_{1},\ldots,X_{n}],\end{equation}
instead of $\Theta_r$ as in (\ref{Th}), one obtains 
a certificate of nonnegativity on $K_S$. This is because 
when using $\theta_r$, the fact that
the (new) approximating sequence $\{f_{\ep r}\}$ lies in $T_S$, also implies that $f$ is nonnegative on $K_S$. Therefore, one may use this property to 
detect whether some given $f$ is nonnegative on $K_S$.

(iii) Finally, we address the issue of identifying the factors that influence 
the degree $r$ up to which one has to perturb $f$ to obtain an s.o.s.
We find that $r$ depends only on $\ep$, the dimension $n$, the degree and the size of the coefficients of $f$, but {\it not} on the explicit choice of $f$.

\textit{Link with related results}. The s.o.s. approximation $f+\ep\Th_{r}$ in (\ref{Th}) resembles the one in (\ref{smalltheta}) recently introduced by the first author in \cite{las3}, for polynomials nonnegative on the \textit{whole} $\R^{n}$; 
with $\th_{r}$ instead of $\Th_{r}$, 
it is proven in \cite{las3} that given a globally nonnegative polynomial $f$ and $\ep>0$, the polynomial $f+\ep\th_{r}$ is s.o.s. provided $r$ is large enough (and we also have $\Vert f+\ep\th_{r} -f\Vert_{1}\rightarrow 0$ as $\ep\rightarrow 0$). Notice that this latter result is also a certificate of nonnegativity on $\R^{n}$ and is more than a denseness result for the $l_{1}$-norm. Indeed, it also shows that every nonnegative polynomial can be approximated by s.o.s. polynomials \textit{uniformly on compact sets}, a nice additional property.

So a polynomial $f$ nonnegative on $\R^{n}$ (hence also on $[-1,1]^{n}$) could be approximated either by $f_{\ep r}=f+\ep\Th_{r}$ or by $f_{\ep r}=f+\th_{r}$ for sufficiently large $r\in\N$; in both cases $\Vert f-f_{\ep r}\Vert_{1}\rightarrow 0$ as $\ep \rightarrow 0$. However, the former approximation is \textit{not} a certificate of nonnegativity of $f$; in particular, it looses the nice property of uniform approximation on compact sets possessed by the latter.

In other words, the s.o.s. approximation $f+\ep\Th_{r}$ is indeed specific for polynomials nonnegative on $[-1,1]^{n}$. For polynomials nonnegative on $\R^{n}$, the s.o.s. approximation $f+\ep\th_{r}$ (although a little more complicated than $f+\ep\Th_{r}$) should be prefered.

The above mentioned Moment Problem for a finite set of polynomials $S\subset\R[X_1,\ldots,X_n]$ is discussed in e.g. \cite{kuh,kuh2}, where the authors ask wether for each polynomial $f$ nonnegative on the corresponding basic closed semi-algebraic 
set $K_S$, there exists some polynomial $q\in\R[X_1,\ldots,X_n]$ such 
that for every $\ep >0$, the polynomial $f+\ep q$ lies in the preordering $T_S$ generated by $S$. This is still an open problem. Our result is weaker, as the polynomial $q$ ($=\Theta_r$ or $\theta_r$) depends on $\ep$ via its degree $r$.

Finally, the degree bounds that we discuss here have been already investigated
in \cite{net} in a similar context, but for the approximations obtained in \cite{las3}.

The paper is organized as follows. After introducing some notation and definitions
in \S \ref{notation}, our results are presented in \S \ref{main1} for s.o.s. approximations of polynomials nonnegative on $[-1,1]^n$, in \S \ref{main2} for related results on polynomials nonnegative on a basic closed semi-algebraic set $K_S\subset\R^n$, and in \S \ref{main3} for results on the degree bounds. For ease of exposition, some technical proofs have been postponed in an Appendix in \S \ref{appendix}.

\section{Notations and definitions}
\label{notation}

Let $\R[X]:=\R[X_{1},\ldots,X_{n}]$ denote the ring of real polynomials, $\A_{r}$ the finite dimensional subspace of polynomials of degree at most $r$ and $s(r)={n+r \choose n}$ its dimension. Let $\A_{r}^{\sos}\subset \A_{r}$ be the space of s.o.s. polynomials of degree as most $r$.

We always fix the canonical monomial basis for $\A_{r}$ and $\R[X_{1},\ldots,X_{n}]$, if we consider them as real vector spaces. For $\al\in\N^{n}$, we write $X^{\al}$ for $X_{1}^{\al_{1}}\cdots X_{n}^{\al_{n}}$, and $|\alpha|$ for $\sum_{i=1}^{n}\alpha_{i}$.

A linear form $L$ on $\R[X]$ is said to have a \textit{representing measure} $\mu$ if $$L(f)=\int_{\R^{n}} f d\mu \qquad \forall f\in\R[X].$$ This is the same as saying that the sequence of values of $L$ on the canonical monomial basis is the moment sequence of this measure $\mu$.

Of course not every linear form has a representing measure. However, there is a \textit{sufficient} condition to ensure that it is indeed the case.

\begin{definition}\label{abval} A function $\ph\colon\N^{n}\rightarrow\R_{+}$ is called an \textit{absolute value} if  \begin{itemize}
    \item[(i)] \( \varphi (0)=1; \)
    \item[(ii)] \( \varphi (\alpha + \beta) \leq \varphi (\alpha) \varphi
    (\beta)\)\ for all \(\alpha , \beta \in\N^{n}\).
    \end{itemize}
\end{definition}

The following result is stated in Berg et al. \cite{bcr}.

\begin{theorem}\label{mom}
Let $L$ be a linear form on $\R[X]$ such that $L(p^{2})\geq 0$ for all $p\in\R[X]$. If there is an absolute value $\ph$ and a constant $C>0$ such that $| L(X^{\alpha}) | \leq C\ph (\alpha)$ for all $\alpha\in\N^{n}$, then $L$ has exactly one representing measure $\mu$ on $\R^{n}$. The support of $\mu$ is contained in the set 
$\left\{x\in\R^{n}\mid\: \vert x^{\alpha}\vert\leq \ph(\alpha)\ \forall\alpha\in\N^{n}\right\}$.
\end{theorem}

For a finite set $S=\left\{g_{1},\ldots,g_{s}\right\}$ of polynomials, denote by $K_{S}$ the basic closed semi-algebraic set $K_{S}:=\left\{x\in\R^{n}\mid g_{i}(x)\geq 0 \ i=1,\ldots,s\right\}$, and by $T_{S}$ the preordering generated by $S$, i.e the set of all finite sums of polynomials of the form $$\si_{e}g_{1}^{e_{1}}\cdots g_{s}^{e_{s}},$$ where $e\in\left\{0,1\right\}^{s}$ and $\si_{e}$ is s.o.s. Further, let $T_{r}$ be the set of all finite sums of such elements  $\si_{e}g_{1}^{e_{1}}\cdots g_{s}^{e_{s}}$ of degree at most $r$. Note that this is different from $T_{S}\cap \A_{r}$ in general, as cancellation of leading forms could result in a polynomial of degree at most $r$, without the single polynomials having this property.

For the degree bound issue addressed in \S \ref{main3}, one needs some elementary notions from the theory of real closed fields and valuation theory. Given a real closed extension field \(R\) of $\R$, denote  by \(\OO\) the convex hull of \(\Z\) in \(R\),  i.e 
\[\OO = \{ x\in R \mid \exists m \in \N :\ |x| \leq m \}.\]
  \(\OO\) is a valuation ring of \(R\) with maximal
   ideal \[ \M = \{x\in\OO \mid \forall n \in \N \setminus \{0\}: \ |x| \leq
   \frac1n \}.\] Let \(\overline{R}:= \OO \slash \M \) denote the residue
   field  and \(\sigma: \OO \to \overline{R}\) the order preserving residue map.
   We have \(\overline{R}=\R\) and \(\sigma\) is the identity on \(\R\). In fact, for every
   \(\beta\in\OO\) there is exactly one \(b\in\R\) such that
   \(\beta\equiv b \mod \M\).

\section{Main results.}

In this section we prove our main results, whereas for ease of exposition, some technical proofs are postponed in \S \ref{appendix}. We first consider polynomials 
nonnegative on the unit ball $[-1,1]^n$.

\subsection{Nonnegativity on the ball $[-1,1]^{n}$}
\label{main1}
We begin with the following result of its own interest.
\begin{theorem}\label{th1}Let $f\in\R[X]$ be a polynomial of degree $r_{f}$, and let $\Th_{r}\in\R[X]$ be as in $(\ref{Th})$. Let $r_{f}\leq 2r\in\N$ be fixed and consider the semidefinite program \begin{equation}\label{sem1}\min_{L}\left\{L(f)\mid L:\mathcal{A}_{2r}\rightarrow \R \mbox{ linear, }L(\Th_{r})\leq 1, \ L(h^{2})\geq 0\  \forall h\in\mathcal{A}_{r} \right\}=:\ep_{r}^{*}.\end{equation}
Then \begin{itemize} \item[(i)] $\ep_{r}^{*} <\infty$ and $(\ref{sem1})$ is solvable, i.e. $\ep_{r}^{*}=L(f)$ for some feasible $L$.
\item[(ii)] The polynomial $f_{\ep r}:=f+\ep\Th_{r}$ is s.o.s. if and only if $\ep\geq -\ep_{r}^{*}$.
\end{itemize}
\end{theorem}
(Note that the condition $L(h^{2})\geq 0\  \forall h\in\mathcal{A}_{r}$ translates to the positive semidefiniteness of the matrix which represents the bilinear form $(p,q)\mapsto L(pq)$. Therefore $(\ref{sem1})$ is an SDP.)

\begin{proof}
(i) The set of feasible solutions for $(\ref{sem1})$ is nonempty, take for example the zero form. So $\ep_{r}^{*}<\infty$. Furthermore, the set of feasible solutions is compact (if we consider each linear form on $\A_{2r}$ as the s(2r)-vector of its values on the monomial basis). Indeed, the constraint $L(\Th_{r})\leq 1$ implies that $$L(1)\leq 1; \quad L(X_{i}^{2r})\leq 1, \quad i=1,\ldots,n.$$ As $L(p^{2})\geq 0$ for all $p\in\A_{r}$, by Lemma $\ref{lem1}$ and Lemma $\ref{lem3}$ from the appendix, one has $|L(X^{\alpha})|\leq 1$ for all $|\alpha|\leq 2r$. So the set of feasible solutions in $\R^{s(2r)}$  is bounded. As it is obviously closed as well, it is compact. Since the objective function is linear and therefore continuous, there always exists an optimal solution.

(ii) By definition, the minimum value $\ep_{r}'$ for which $f_{\ep r}$ is s.o.s. is given by \begin{equation}\label{sem2} \ep_{r}'=\min_{\ep}\left\{\ep\mid f+\ep\Th_{r}\in \A_{2r}^{\sos}\right\}.\end{equation} But $(\ref{sem2})$ is an SDP whose dual reads $$\max_{L}\left\{L(f)\mid L:\mathcal{A}_{2r}\rightarrow \R \mbox{ linear, }-L(\Th_{r})\leq 1, \ L(h^{2})\leq 0\  \forall h\in\mathcal{A}_{r} \right\}.$$ Equivalently, with the change of variable $L\rightarrow -L$, \begin{equation}\label{sem3}-\min_{L}\left\{L(f)\mid L:\mathcal{A}_{2r}\rightarrow \R \mbox{ linear, }L(\Th_{r})\leq 1, \ L(h^{2})\geq 0\  \forall h\in\mathcal{A}_{r} \right\}.\end{equation}
One next proves that there is no \textit{duality gap} between the respective primal and dual problems $(\ref{sem2})$ and $(\ref{sem3})$, that is, their respective optimal values are equal.

Let $\mu$ be a measure on $\R^{n}$ with all moments up to order $2r$ finite and with a strictly positive density. One may scale $\mu$ to satisfy $\int_{\R^{n}}\Th_{r}d\mu < 1$. Let $L$ be integration with respect to $\mu$. As $\mu$ has strictly positive density, we must have $L(p^{2})> 0$ for all $p\in\A_{r}\setminus\{0\}$, and so $L$ is a strictly feasible solution for the SDP in $(\ref{sem3})$, that is, Slater's condition holds, which in turn implies that both SDP problems in $(\ref{sem2})$ and $(\ref{sem3})$ have the same optimal value $\ep_{r}'=-\ep_{r}^{*}$; see e.g. \cite{van}.

So the \textit{only if} part in (ii) follows from the definition of $\ep_{r}'$. Now let $\ep\geq-\ep_{r}^{*}$ and write
$$f+\ep\Th_{r}=f-\ep_{r}^{*}\Th_{r}+(\ep+\ep_{r}^{*})\Th_{r},$$ and use that $f-\ep_{r}^{*}\Th_{r}$ as well as $(\ep+\ep_{r}^{*})\Th_{r}$ are s.o.s. to obtain the result.
\end{proof}

Observe that $\ep_{r}^{*}=0$ whenever $f$ is a s.o.s., because then $L(f)\geq 0$ for every feasible $L$ and the zero linear form is feasible. If $f$ is not s.o.s. (so $\ep_{r}^{*}<0$), then the inequality constraint $L(\Th_{r})\leq 1$ in $(\ref{sem1})$ can be replaced with the equality constraint $L(\Th_{r})=1$, since by linearity, given a feasible solution $L$ with $L(\Th_{r})<1$ and with value $L(f)<0$, one always obtains a better feasible solution $L'=\rh L$ with $L'(\Th_{r})=1$ (note that $L(\Th_{r})=0$ implies $L=0$).

Next, we obtain the following crucial result.
\begin{theorem}\label{th2}
Let $f\in\R[X]$ be a polynomial of degree $r_{f}$, nonnegative on $[-1,1]^{n}$, and let $\Th_{r}\in\R[X]$ be as in $(\ref{Th})$. Let $\ep_{r}^{*}$ be the optimal value of the semidefinite program defined in $(\ref{sem1})$, for all $2r\geq r_{f}$. Then $\ep_{r}^{*}\rightarrow 0$ as $r\rightarrow \infty$.
\end{theorem}
\begin{proof}
 From Theorem $\ref{th1}$, $\ep_{r}^{*}=L^{(r)}(f)\leq 0$ for some optimal solution $L^{(r)}$ of the semidefinite program $(\ref{sem1})$, whenever $2r\geq r_{f}$. From the proof of Theorem $\ref{th1}$, it follows that $|L^{(r)}(X^{\alpha})|\leq 1$ for all $\alpha\in\N^{n}$ with $|\alpha|\leq 2r$. Next, complete the vector of the values of $L^{(r)}$ on the monomial basis of $\A_{r}$ with zeros to make it an element in $\R^{\N^{n}}$, and in fact even an element of $[-1,1]^{\N^{n}}$. By Tychonoff's Theorem, we find a subsequence $r_{k}$ such that the sequence $L^{(r_{k})}$ converges to some ${\bf y}^{*} \in [-1,1]^{\N^{n}}$ in the product topology, and in particular pointwise convergence holds, i.e. \begin{equation}\label{con} L^{(r_{k})}(X^{\alpha})\rightarrow {\bf y}^{*}_{\alpha} \qquad \forall \alpha\in\N^{n}.\end{equation}

 Let $L^{*}$ be the linear form on $\R[X]$ defined by $L^{*}(X^{\alpha}):={\bf y}^{*}_{\alpha}$. From the pointwise convergence in $(\ref{con})$ we obtain $L^{*}(p^{2})\geq 0$ for all $p\in\R[X]$. This, together with ${\bf y}^{*}\in [-1,1]^{\N^{n}}$, implies that $L^{*}$ has a representing measure $\mu^{*}$ with support contained in $[-1,1]^{n}$ (see Theorem $\ref{mom}$). Now again from the pointwise convergence $(\ref{con})$,
 $$ \ep_{r_{k}}^{*}=L^{(r_{k})}(f)\rightarrow L^{*}(f)=\int_{[-1,1]^{n}}fd\mu^{*}\geq 0,$$ where the inequality uses nonnegativity of $f$ on $[-1,1]^{n}$. Since all $\ep_{r}^{*}\leq 0$, we get $\ep_{r_{k}}^{*}\rightarrow 0$. And as the converging subsequence $r_{k}$ was arbitrary, this shows the desired result.
\end{proof}

Therefore, we finally obtain:
\begin{corollary}\label{app}
Let $f\in\R[X]$ be a polynomial nonnegative on $[-1,1]^{n}$ and
let $\Th_{r}\in\R[X]$ be as in $(\ref{Th})$. Let $\ep>0$ be fixed. Then there exists some  $r(f,\ep)\in\N$ such that for
every $r\geq r(f,\ep)$, the polynomial $f_{\ep r}:=f+\ep\Th_{r}$
is a s.o.s.
\end{corollary}
\begin{proof}
From Theorem $\ref{th2}$ we know that the sequence $\{\ep_{r}^{*}\}$ with $\ep_{r}^{*}$ defined in $(\ref{sem1})$ converges to $0$ as $r\rightarrow\infty$. So there is an $r(f,\ep)$ such that for all $r\geq r(f,\ep)$ we have $\ep_{r}^{*}\geq -\ep$. By Theorem $\ref{th1}$ the polynomial $f-\ep_{r}^{*}\Th_{r}$ is a s.o.s., and so $$f+\ep\Th_{r}=f-\ep_{r}^{*}\Th_{r} + (\ep+\ep_{r}^{*})\Th_{r}$$ is a s.o.s. as well, since  $(\ep+\ep_{r}^{*})\Th_{r}$ is also a s.o.s. ($\ep_{r}^{*}\geq-\ep$).
\end{proof}

Corollary $\ref{app}$ refines the \textit{denseness} result of
Berg \cite{ber},  because it provides an \textit{explicit}
approximation sequence. In addition, this approximation sequence
is extremely simple, as the perturbation polynomial $\Th_{r}$
contains only the constant and the marginal monomials
$X_{i}^{2r},\ i=1,\ldots,n.$ In addition, it provides a
\textit{certificate} of nonnegativity of $f$ on $[-1,1]^{n}$;
indeed, if $x\in [-1,1]^{n}$, then for every $r\geq r(f,\ep)$ one
has $f(x)+\ep\Th_{r}(x)\geq 0$. Letting $\ep\rightarrow 0$ yields
$f(x)\geq 0$.

It is straightforward to extend Corollary $\ref{app}$ to the case
of a polynomial $f$  nonnegative on the ball
$[-l,l]^{n}\subset\R^{n}$ for some $l>0$. Indeed, it suffices to
apply Corollary $\ref{app}$ to the polynomial $x\mapsto
g(x):=f(lx)$ which is nonnegative on $[-1,1]^{n}$. In this case
the polynomial $f+\ep(1+\sum_{j=1}^{n}(X_{j}/l)^{2r})$ provides an
s.o.s. approximation.
\bigskip

In some specific examples, one may even obtain a more precise
result. Namely, given $r$ \textit{fixed}, one may provide an
explicit bound $\ep_{r}>0$, such that the polynomial $f_{\ep
r}:=f+\ep\Th_{r}$ is s.o.s. This is illustrated in the following
nice two examples, kindly provided by Bruce Reznick.

\begin{example}\label{ex1}
Consider the univariate polynomial $f=1-X^{2}$, obviously nonnegative on $[-1,1]$. If $\ep\geq \ep_{r}^{*}:= (r-1)^{r-1}/r^{r}$, the polynomial $$f_{\ep r}:=1-X^{2}+\ep X^{2r}$$ is globally nonnegative and therefore a s.o.s. Indeed, its minimum occurs when $-2x+2r\ep x^{2r-1}=0$, i.e. at $x_{r}:=(1/r\ep)^{1/(2r-2)}$. Hence, the value at $x_{r}$ is $$1-x_{r}^{2}+\ep x_{r}^{2}x_{r}^{2r-2}=1-x_{r}^{2}(r-1)/r,$$ which is nonnegative if and only if $$x_{r}^{2}\leq r/(r-1) \Leftrightarrow x_{r}^{2r-2}\leq (r/(r-1))^{r-1} \Leftrightarrow 1/(r\ep) \leq (r/(r-1))^{r-1},$$ i.e. if and only if $\ep\geq (r-1)^{r-1}/r^{r}=\ep_{r}^{*}$.
\end{example}
\begin{example}\label{ex2}
On the other hand, consider the Motzkin polynomial $f=1+X^{2}Y^{2}(X^{2}+Y^{2}-3)\in\R[X,Y]$ which is nonnegative but \textit{not} a s.o.s. Then, for all $r\geq 3$ and $\ep:=2^{4-2r}$, the polynomial $f_{\ep r}:=f+\ep X^{2r}$ is a s.o.s., and $\Vert f-f_{\ep r}\Vert_{1}\rightarrow 0$ as $r\rightarrow\infty$. To prove this, write $$f=(XY^{2}+X^{3}/2-3X/2)^{2}+p,$$ where $p=1-(X^{3}/2-3X/2)^{2}=(1-X^{2})^{2}(1-X^{2}/4).$ Next, the univariate polynomial $q=p+2^{4-2r}X^{2r}$ is nonnegative on $\R$, hence a sum of squares. Indeed, if $x^{2}\leq 4$, then $p\geq 0$ and so $q\geq 0$. If $x^{2}>4$ then $|p(x)|\leq (x^{2})^{2}x^{2}/4=x^{6}/4.$ From $$q(x)\geq 2^{4-2r}x^{2r}-|p(x)|\geq \frac{x^{6}}{4}((x^{2}/4)^{r-3}-1),$$ and the fact that $n\geq 3, x^{2}>4$, we deduce that $q(x)\geq 0$.
\end{example}

In Example $\ref{ex1}$, one approximates $1-X^{2}$ (uniformly on $[-1,1]$) by the s.o.s. $1-X^{2}+\ep X^{2r}$. In Example $\ref{ex2}$, the Motzkin polynomial can also be approximated in the $l_{1}$-norm by $f+\ep(X^{2r}+Y^{2r})$, but not uniformly on compact sets. For the latter property to hold, one needs the perturbation $f+\ep\sum_{j=1}^{n}\sum_{k=0}^{r}X_{i}^{2k}/k!$ introduced in \cite{las3}.
\bigskip

\subsection{Nonnegativity on basic closed semi-algebraic sets}
\label{main2}
We next prove the second announced result,
namely the approximation of polynomials nonnegative on basic closed
semi-algebraic sets. Let $S\subset\R[X]$ be a finite set of
polynomials and suppose the Moment Problem is solvable for $S$,
which means that every linear form on $\R[X]$ which is nonnegative
on the preordering $T_{S}$, is integration with respect to some
measure on $K_{S}$. Further suppose $K_{S}$ has nonempty interior, and
let $f\in\R[X]$ be nonnegative on $K_{S}$.

With same notation as in \S \ref{main1}, 
consider the semidefinite program \begin{equation}\label{sempre1}
\ep_{r}^{*}:=\min_{L}\left\{L(f)\mid L\colon\A_{2r}\rightarrow\R \ {\rm linear}, \ L(\Th_{r})\leq 1,\ L(t)\geq 0 \ \forall t\in T_{2r}\right\}.
\end{equation}
Its dual reads \begin{equation}\label{sempre2} \max_{\ep}\left\{\ep\mid f-\ep\Th_{r}\in T_{2r}\right\}.\end{equation} 

Proceding exactly as in the proof of Theorem \ref{th1}, one constructs a strictly feasible solution for $(\ref{sempre1})$ as integration with respect to some (suitably scaled) measure on a ball in $K_{S}$. Hence, with same arguments, 
the SDP (\ref{sempre1}) is also always 
solvable (note that $\A_{2r}^{\sos}\subseteq T_{2r}$),
and there is no duality gap between the SDPs (\ref{sempre1}) and (\ref{sempre2}), i.e., their optimal values are equal.

Again, every sequence of optimal solutions for $(\ref{sempre1})$ (with $r$ growing) has a subsequence that converges pointwise to some ${\bf y}^{*}\in [-1,1]^{\N^{n}}$ wich is the moment sequence of some measure on $K_{S}$, this time using the fact that the moment problem holds for $S$. So, as in the proof of Theorem \ref{th2}, 
the sequence $\{\ep_{r}^{*}\}$ converges to $0$, since $f$ is nonnegative on $K_{S}$. Hence, as in Corollary $\ref{app}$, we get the following result:

\begin{corollary}\label{app2}
Let $S\subset\R[X]$ be a finite set of polynomials and suppose that the Moment Problem is solvable for $S$. Further, suppose that $K_{S}$ has a nonempty interior. Let $f\in\R[X]$ be nonnegative on $K_{S}$ and let $\Th_{r}\in\R[X]$ be as in $(\ref{Th})$. Let $\ep>0$ be fixed. Then there is some $r(f,\ep,S)$ such that for every $r\geq r(f,\ep,S)$, the polynomial $f_{\ep r}:=f+\ep\Th_{r}$ lies in $T_{S}$.
\end{corollary}
Note that the pointwise limit ${\bf y}^{*}$ from above is the moment sequence of a measure on $K_{S}$ as the Moment Problem holds for $S$, but on the other hand it is also the moment sequence of a measure on $[-1,1]^{n}$, as ${\bf y}^{*}\in [-1,1]^{\N^{n}}$ (Theorem $\ref{mom}$). But by Theorem  $\ref{mom}$, ${\bf y}^{*}$ is the moment sequence of exactly one measure. So the measure must be supported by $K_{S}\cap [-1,1]^{n}$. This leads to the fact that in Corollary $\ref{app2}$, the polynomial $f$  must only be nonnegative on $K_{S}\cap [-1,1]^{n}$ for the statement to hold. So for example if $[-1,1]^{n}\cap K_{S} =\emptyset$, it holds for every polynomial $f$.

However, notice that "$f+\epsilon\Th_r$ {\it lies in $T_S$}" provides a 
certificate of nonnegativity of $f$ on $K_S\cap [-1,1]^n$ only, and {\it not} on $K_S$. So Corollary \ref{app2} is useful when one already knows that $f$ is nonnegative on $K_S$ and one wishes to obtain an $l_1$-norm approximation in $T_S$. If one wishes to test whether $f$ is indeed nonnegative on $K_S$, then the following result provides a certificate of nonnegativity on $K_S$.

\begin{corollary}\label{app3}
Let $S\subset\R[X]$ be a finite set of polynomials and suppose that the Moment Problem is solvable for $S$. Further, suppose that $K_{S}$ has a nonempty interior. Let $f\in\R[X]$ be nonnegative on $K_{S}$ and let $\th_{r}\in\R[X]$ be as in $(\ref{smalltheta})$. Let $\ep>0$ be fixed. Then there is some $r(f,\ep,S)$ such that for every $r\geq r(f,\ep,S)$, the polynomial $f_{\ep r}:=f+\ep\th_{r}$ lies in $T_{S}$.
\end{corollary}
The proof is similar to that of Corollary \ref{app2}, except that in the semidefinite 
program (\ref{sempre1}) we now have the constraint  
$L(\theta_r)\leq 1$ (instead of $L(\Theta_r)\leq1$).
In this case, every sequence of optimal solutions for $(\ref{sempre1})$ (with $r$ growing) has a subsequence that converges pointwise to some ${\bf y}^{*}\in \R^{\N^{n}}$ (rather than ${\bf y}^*\in[-1,1]^{\N^{n}}$). To prove this result, and as one cannot use
Theorem \ref{mom} any more, one now invokes Nussbaum's result
\cite{nuss} on moment sequences, which, in the present context, states that if
\[\sum_{i=1}^n\sum_{k=1}^\infty \,L(X_i^{2k})^{-1/2k}\,=\,+\infty,\qquad i=1,\ldots,n,\]
then $L$ is integration with respect to some measure on $\R^n$; see also Berg \cite[Theorem 8]{ber}. The rest of the proof is identical.\\

That Corollary \ref{app3} provides a certificate of nonnegativity of $f$ on $K_S$, follows from the fact that $\theta_r(x)$ is bounded by $\sum_{i=1}^n\mathrm{exp}(x_i^2)$, for all $x\in\R^n$. Therefore, fix $x\in K_S$; as $f+\epsilon\theta_r$ lies in $T_S$, one has
$f(x)+\epsilon\theta_r(x)\geq0$. Letting $\epsilon\to 0$ yields $f(x)\geq0$, the desired result.\\

The result in Corollary $\ref{app2}$ (resp. in Corollary \ref{app3}) is weaker than the condition $f+\ep q\in T_S$ for some fixed $q$ and all $\ep>0$, as our $\Th_{r}$ (resp. $\theta_r$) depends on $\ep$ (via $r$). Whether the Moment Problem implies even this stronger version is an open problem, see for example \cite{kuh,kuh2}.
\bigskip

\subsection{The degree of the perturbation}
\label{main3}
We are now concerned with the last announced result. We prove that  the
degree $r(f,\ep)$ in Corollary $\ref{app}$ does {\it not} depend on the
{\it explicit choice} of the polynomial $f$ but only on 

$\bullet$ $\ep$ and the dimension $n$,

$\bullet$ the degree and the size of the coefficients of $f$.\\

Therefore, if we fix these four parameters, we find an $r$ such
that the statement of Corollary $\ref{app}$ holds for any $f$
nonnegative on $[-1,1]^{n}$, whose degree and size of the
coefficients do not exceed the fixed parameters. 

We first generalize Corollary $\ref{app}$ to real closed extension fields of $\R$ and then use the result in an ultrapower of $\R$. This approach towards degree bounds is similar to the one in \cite{pre}. 

Let $\Th_{r}$ be as in (\ref{Th}). We first write the strict duality of the SDP problems $(\ref{sem2})$ and $(\ref{sem3})$ as a first order logic formula in the language of ordered rings with coefficients from $\R$. We just say that for every polynomial $f$ of some fixed maximum degree $2r$, there is a linear form $L$  on $\A_{2r}$  (indeed a $s(2r)$-tuple of values) which is nonnegative on $\A_{2r}^{\sos}$ and which is less than or equal to $1$ on $\Th_{r}$. We also demand that all the values of $L$ on the monomial basis are bounded by $1$ (as we have seen, this follows from the other conditions anyway).
Further, we say that there exists some $\ep$ such that $f+\ep\Th_{r}$ is a s.o.s and $\ep=-L(f)$ with $L$ from above. All this can be done, using the known fact that every polynomial in $\A_{2r}^{\sos}$ is already a sum of $s(2r)$ squares of polynomials from $\A_{r}$.

So, by Tarski's Transfer Principle, 
for every $r\in\N$, this formula holds in every real closed extension field of $\R$.
We use this in the following theorem:

\begin{theorem}\label{ext}
Let $R$ be a real closed extension field of $\R$, and denote by $\OO$ the convex hull of $\Z$ with respect to the unique ordering in $R$. Let $\M$ denote the unique maximal ideal in the valuation ring $\OO$, and fix some $\ep\in R, \ep>0$ and $\ep \notin\M$. Suppose $f\in\OO[X]$ is nonnegative on $[-1,1]^{n}\subset R^{n}$. Then there exists $r\in\N$ such that  the polynomial $f_{\ep r}=f+\ep\Th_{r}$ is a s.o.s. in $R[X]$.
\end{theorem}
\begin{proof}
 Let $\overline{f}$ be the real polynomial obtained from $f$ by applying the residue map $\si\colon\OO\rightarrow \OO / \M=\R$ to the coefficients of $f$.  As $f\geq 0$ on $[-1,1]^{n}\subset R^{n}$, we have $\overline{f}\geq 0$ on $[-1,1]^{n}\subset\R^{n}$.

 Next, consider the SDP problems from $(\ref{sem1})$ associated with $\overline{f}$.
 From Theorem $\ref{th2}$, there exists some $r$ such that $\ep_{r}^{*}> - \si(\ep)$ ($\ep >0, \ep\notin\M$ implies $\si(\ep)> 0$). With that $r$ fixed, we now use that the formula described above holds in $R$. That is, we first get a linear form $L$ on the subspace of polynomials of $R[X]$ with degree at most $2r$, whose values on the monomial basis are bounded by $1$ (and therefore, are in $\OO$), which is nonnegative on the s.o.s. polynomials. Further, we also have $L(\Th_{r})\leq 1$. In addition, we get an $\ep'$ such that $f+\ep'\Th_{r}$ is a s.o.s. in $R[X]$ and $\ep'=-L(f)$.

 But now, we can apply the residue map $\si$ to the values of $L$ on the monomial basis and get a linear form $\overline{L}$ wich is feasible for the optimization problem from $(\ref{sem1})$ associated with $\overline{f}$ and $r$. So $$-\si(\ep)< \ep_{r}^{*}\leq \overline{L}(\overline{f})=\si({L(f)})=-\si({\ep'}).$$ This shows $\ep'<\ep$, and as $f+\ep'\Th_{r}$ is a s.o.s. in $R[X]$, so is $f+\ep\Th_{r}$.
\end{proof}

Once we have this result, the rest follows from a  standard ultrapower
argument. We use the result in $$\R^{*}= (\prod_{\N} \R )/
\mathcal{U},$$ where $\mathcal{U}$ is a non-principal ultrafilter
on $\N$. 

Fix some $\ep\in\R$, $\ep>0$, and define
by a first order logic formula $\Ph$ in the language of ordered rings,
the set of all polynomials $f$ of degree at most $d$, with coefficients
bounded by some $N\in\N$, and which are nonnegative on $[-1,1]^{n}$.

Next, for every $r\in\N$, define by a formula $\ph_{r}$, the set of all
polynomials $f$ of degree at most $d$, such that $f+\ep\Th_{r}$ is
a s.o.s. 

Notice that boundedness of the
coefficients of a polynomial $f$ by some $N\in\N$, implies
$f\in\OO[X]$, and so, by Theorem $\ref{ext}$, one has
$$\Ph \rightarrow \bigvee_{r\in\N} \ph_{r}.$$ Now the
$\aleph_{1}$-saturation of $\R^{*}$ yields $$\Ph\rightarrow
\ph_{r'}$$ for some $r'$ depending on the formulas used, i.e. on
$d,N,n,\ep$. Therefore, in $\R^{*}$ one may choose the degree $r$
in Theorem $\ref{ext}$ to depend only on $d,N,n,\ep$. As this can
be again formulated as a first order logic formula, it holds in $\R$ as well:
\begin{theorem}
Let $n,N,d\in\N$ and $\ep\in\R_{>0}$ be given.
Then there exists $r=r(n,N,d,\ep)\in\N$ such that for  every
$f\in\R[X_{1},\ldots,X_{n}]$ of degree at most $d$, with
coefficients bounded by $N$, and nonnegative on $[-1,1]^{n}$,
the polynomial $f+\ep\Th_{r}$ is a s.o.s. (and so are $f+\epsilon\Theta_{r'}$
for all $r'\geq r$).
\end{theorem}

\section{appendix}
\label{appendix}
In this section we derive auxiliary results that are helpful in the proofs of the main section.
\begin{lemma}\label{lem1}
Let $n=1$ and let $L\colon\A_{2r}\rightarrow\R$ be a linear form such that $L(p^{2})\geq 0$ for all $p\in\A_{r}$. Then $L(X^{2k})\leq \max [L(1),L(X^{2r})]$ for all $k=0,\ldots,r.$
\end{lemma}
\begin{proof}
The proof is by induction on $r$. Indeed for $r=0$ and $r=1$ the statement is trivial. So we assume the statement of Lemma $\ref{lem1}$ is true for some $r$ and we prove it for $r+1$.

Let $L$ be a linear form on $\A_{2r+2}$ as stipulated. From $L(p^{2})\geq 0$ for all $p\in\A_{r+1}$ we have \begin{equation}\label{ab}
L(X^{2r})^{2} \leq L(X^{2r+2})L(X^{2r-2}).                                                                                 \end{equation}
By the induction hypothesis, we have $$L(X^{2k})\leq \max [L(1),L(X^{2r})],\qquad k=0,\ldots,r.$$ Suppose first that $L(1)=\max [L(1),L(X^{2r})]$. Then obviously \linebreak $L(X^{2k})\leq \max [L(1),L(X^{2r+2})]$ for all $k\leq r+1$ and we are done. Next, suppose $L(X^{2r})=\max [L(1),L(X^{2r})]$. Then from $(\ref{ab})$ we obtain $$L(X^{2r})^{2}\leq L(X^{2r+2})L(X^{2r-2})\leq L(X^{2r+2})L(X^{2r}),$$ so that $L(X^{2r})\leq L(X^{2r+2})$. Therefore again $L(X^{2k})\leq \max [L(1),L(X^{2r+2})]$ for all $k=0,\ldots,r+1,$ the desired result.
\end{proof}

\begin{lemma}\label{lem2} Let $n=2$ and \(L:\A_{2r}\rightarrow \R\) be a linear form and suppose \(L(p^{2})\geq 0\) for all \(p\in\A_{r}\). Then all values \(L(X^{2\alpha})\) where \(0\leq |\alpha | \leq r \) are bounded by \(\max_{k=0,...,r} \max \{L(X_{1}^{2k}),L(X_{2}^{2k})\}.\)
   \end{lemma}
   \begin{proof} What we will actually show is that all \(L(X^{2\alpha})\),
   where \(|\alpha | = k\), are bounded by \(\max\{ L(X_{1}^{2k}), L(X_{2}^{2k})\}\).

   Let \(p\in\N\) be such that either \(k=2p\) (if \(k\)
   is even) or \(k=2p+1\) (if \(k\) is odd) and define \(\Gamma:=\{\ (2a,2b) \ | \ a + b = k;\ a,b \neq 0\}\). One
   has \(\Gamma = \Gamma_{1}  \cup \Gamma_{2}\) where \begin{eqnarray*}\Gamma_{1}
   & := & \{ \ (k,0)+(k-2i,2i) \ | \ i=1,...,p\ \} \\
   \Gamma_{2}&:=&\{ \ (0,k)+(2j,k-2j) \ | \ j=1,...,p\}.\end{eqnarray*}
   If \(k\) is odd, then this union is disjoint, else \(\Gamma_{1}
   \cap \Gamma_{2} = \{(2p,2p)\}\).
    For \(s:= \max\{ L(X^{\gamma})  \mid \gamma\in \Gamma \ \}\), we get
   \(s=L(X^{\gamma^{*}})\) for some \(\gamma^{*}\in\Gamma_{1} \mbox{ or } \gamma^{*}\in\Gamma_{2}\).

   From $L(p^{2})\geq 0$ for all $p\in\A_{r}$ we have \[L(X^{\alpha + \beta})^{2} \leq L(X^{2\alpha})L(X^{2\beta}).\]

   So in our case, we  obtain \begin{eqnarray}\label{6} L\left(X_{1}^{2k}\right)\cdot L\left(X_{1}^{2k-4i}X_{2}^{4i}\right) &\geq&
   L\left(X_{1}^{2k-2i}X_{2}^{2i}\right)^{2}, \ i=1,...,p \\ \label{7} L\left(X_{2}^{2k}\right) \cdot L\left(X_{1}^{4j}X_{2}^{2k-4j}\right) &\geq&
   L\left(X_{1}^{2j}X_{2}^{2k-2j}\right)^{2}, \ j=1,...,p.\end{eqnarray}

   With \(s_{k}:= \max\{L\left(X_{1}^{2k}\right) , L\left(X_{2}^{2k}\right)\}\), by (\(\ref{6}\)) and (\(\ref{7}\)), one
   gets either \[s_{k}\cdot s\geq L\left(X_{1}^{2k}\right) \cdot L\left(X^{\gamma^{*}}\right) \geq
   L\left(X^{\gamma^{*}}\right)^{2}=s^{2} \]or \[s_{k} \cdot s\geq
   L\left(X_{2}^{2k}\right) \cdot L\left(X^{\gamma^{*}}\right) \geq L\left(X^{\gamma^{*}}\right)^{2} = s^{2}.\] In
   any case \(s_{k} \geq s\). \end{proof}

    \begin{lemma} \label{lem3}
   Let $n$ be arbitrary and \(L:\A_{2r} \rightarrow \R\) be a linear form and suppose \(L(p^{2})\geq 0\) for all \(p\in\A\).
   Assume that for all i=1,...,n and k=0,...,r, the values \(L(X_{i}^{2k})\) are bounded by some \(\tau\). Then all values \(L(X^{\alpha})\), where \(|\alpha|\leq 2r\), satisfy \(|L(X^{\alpha})|\leq \tau.\)

   \end{lemma}
   \begin{proof} We only need to show that all values
   \(L\left(X^{2\alpha}\right)\), where \(|\alpha|\leq r\), are bounded by \(\tau\). Indeed, from $L(p^{2})\geq 0$ for all $p\in\A_{r}$ we have \(L(X^{\alpha + \beta})^{2} \leq L(X^{2\alpha})L(X^{2\beta}),\) and therefore, if all the values \(L\left(X^{2\gamma}\right)\) are bounded by \(\tau\), one gets \( |L\left(X^{\alpha}\right)|\leq \tau \) for all \(0 \leq|\alpha|\leq 2r\).

   The proof is by induction on the number \(n\) of variables.

   \(n=1:\) Nothing is to be shown in this case, as all the values \(L\left(X^{2\alpha}\right)\) are bounded by \(\tau\) by the assumption.

   \(n=2:\) This is an immediate result of Lemma \(\ref{lem2}\).

   \(n-1\rightsquigarrow n, n>2:\) By the induction hypothesis, the
   claim is true for all \(L\left(X^{2\alpha}\right)\), where \(|\alpha| \leq r\) and some \(\alpha_{i}
   =0\). Indeed, \(L\) restricts to a linear form on the ring of polynomials with \(n-1\) indeterminates and satisfies all the assumptions needed. So the induction hypothesis gives the boundedness of all those values \(L\left(X^{2\alpha}\right)\).

   Now take \(L\left(X^{2\alpha}\right)\), where \(|\alpha|\leq r \) and all \(\alpha_{i} \geq 1\). With no loss of generality, assume \(\alpha_{1}\leq \alpha_{2} \leq...\leq\alpha_{n}.\)
   Consider the two elements \[\gamma:=(2\alpha_{1},0,\alpha_{3}+\alpha_{2}-\alpha_{1},\alpha_{4},...,\alpha_{n})\in\N^{n} \mbox { and }\]
   \[\gamma^{'}:=(0,2\alpha_{2},\alpha_{3}+\alpha_{1}-\alpha_{2},\alpha_{4},...,\alpha_{n})\in\N^{n}.\]
   We have \(|\gamma|,|\gamma^{'}|\leq r\) and
   \(\gamma_{2}=\gamma^{'}_{1}=0\). Therefore, by the above
   result, we get \[L\left(X^{2\gamma}\right)\leq\tau \mbox{ and }
   L\left(X^{2\gamma^{'}}\right) \leq\tau.\] As \(L(p^{2})\geq 0\) for all $p\in\A_{r}$  one has \[L\left(X^{2\alpha}\right)^{2}= L\left(X^{\gamma +\gamma^{'}}\right)^{2}\leq
   L\left(X^{2\gamma}\right)\cdot L\left(X^{2\gamma^{'}}\right) \leq\tau^{2},\] which yields
   \[|L\left(X^{2\alpha}\right)|\leq\tau.\] \end{proof}
   
 \section*{Acknowledgements}
Both authors wish to thank M. Schweighofer for many interesting and helpful discussions on the topic. The work of the first author is partly supported by ANR Grant ${\rm NT}\,05-3-41612$, while that of the second author is supported by the Land Baden-W\"urttemberg through a Landesgraduiertenstipendium.

\end{document}